\documentclass[12pt,a4]{article}
\usepackage{graphicx,amssymb,amsmath,amsthm,cite}

\def\be{\begin{equation}}
\def\ee{\end{equation}}
\def\ben{\begin{eqnarray}}
\def\een{\end{eqnarray}}
\def\inoj{i=1 \atop i \ne j}

\def\bmth#1{{\mbox{\boldmath$#1$}}}
\def\ei{{\bmth{e}}_i}
\newcommand{\til}{\tilde}
\newcommand{\ut}{\til{u}}
\newcommand{\utk}{\til{u}^k}
\newcommand{\utkp}{\til{u}^{k+1}}
\newcommand{\utik}{\ut_i^k}
\newcommand{\utikp}{\ut_i^{k+1}}

\newcommand{\la}{\langle}
\newcommand{\ra}{\rangle}

\DeclareMathOperator{\Span}{span}

\newtheorem{lemma}{Lemma}

\newtheorem{proposition}{Proposition}
\newtheorem{property}{Property}
\newtheorem{remark}{Remark}

\def\sumik{\sum_{i=1}^{k}}
\def\sumikp{\sum_{i=1}^{k+1}}

\def\C{\mathbb{C}}
\def\emptyy{\{0\}}
\def\h{\cal{H}}
\def\hk{{\cal{H}}_k}
\def\op{\hat{P}}

\def\kj{k \setminus j}

\def\SV{{\cal{V}}}
\def\SVK{{\cal{V}}_k}
\def\SVKP{{\cal{V}}_{k+1}}

\def\SVKM{{\cal{V}}_{\kj}}
\def\SW{{\cal{W}}}
\def\SWK{{\cal{W}}_k}
\def\SWC{{\cal{W}^\bot}}
\def\SVC{{\cal{V}^\bot}}
\def\EVW{\hat{E}_{\SV \SWC}}
\def\EVV{\hat{E}_{\SV \SVC}}
\def\EVKW{\hat{E}_{\SV_k \SWC}}
\def\EVKWP{\hat{E}_{\SV_{k+1}\SWC}}
\def\EVKWM{\hat{E}_{\SV_{\kj}\SWC}}
\def\utikj{\ut_i^{\kj}}

\def\oo{{\cal{O}}}
\def\nuo{{\cal{N}}(\oo)}
\def\rao{{\cal{R}}(\oo)}

\title{Constructive updating/downdating of oblique projectors:
  a generalization of the Gram-Schmidt process}
\author{Laura Rebollo-Neira\\
Aston University, Birmingham B4 7ET, UK}
\date{}
\begin{document}
\maketitle
\baselineskip = 1\baselineskip
\begin{abstract}
A generalization of the Gram-Schmidt procedure is achieved by 
providing equations for updating and downdating oblique projectors.
The work is motivated by 
the problem of adaptive signal representation outside the 
orthogonal basis setting. The proposed techniques are shown to 
be relevant to the problem of 
discriminating signals produced by 
different phenomena when the order of the signal 
model needs to be adjusted.
\end{abstract}

\section{Introduction}
An operator $\oo$ is a projector if it is idempotent, i.e., if 
it is endowed with property $\oo^2 =\oo$. The projection is 
along (or parallel to) its null space $\nuo$ and onto its range
$\rao$. This entails  that $\oo v = v$ for  
$v \in \rao$ and $\oo v =0$ for $v \in \nuo$. 
If the subspaces $\rao$ and  $\nuo$
are orthogonal the operator is called an
orthogonal projector, which  is the case if and only if 
$\oo$ is self-adjoint. Otherwise it is called oblique projector. 

Oblique projectors, though introduced  early \cite{ob1,ob2}, 
have received less attention than orthogonal projectors.
Nevertheless, quite recently there has been a renewed interest 
in relation to their properties and applications 
\cite{ob3,ald,for1,tak,gro,for2,cen,go1,go2,book}.  
In particular, oblique projectors have 
been shown to be of significant relevance to signal processing 
techniques \cite{sp1,sp3,yon1,ref}.   
The present effort is very much motivated by problems 
arising in the area of signal representation outside
the traditional orthogonal basis setting 
\cite{mal,relo,rem,andre1}.
In such a context a signal $f$, represented mathematically 
as an element of a vector space, is approximated as a linear 
expansion of the form 
\be \label{eq1} f^k=\sum_{i=1}^k c_i v_i. \ee
The vectors $v_i$ in \eqref{eq1} are sometimes sequentially 
fed or chosen according to some optimality
criterion. In such situations one needs to be in a position 
to effectively 
adapt the coefficients of the linear superposition so as to
account for the possibility of changes in the model. 
This may entail a)increasing the order $k$ of the model 
by incorporating new terms in the expansion 
b)reducing the order by eliminating some terms in the 
expansion c)replacing some of the vectors 
in \eqref{eq1} by different ones.

Assuming that the signal space is an inner product space, 
for $f^k$  given in \eqref{eq1} to be the best approximation 
of a signal $f$ in a minimum distance sense, the coefficients 
in \eqref{eq1} should be calculated  in such a way
that $f^k$ is the orthogonal projection of $f$ onto 
$\Span\{v_i\}_{i=1}^k$. This is the main reason for 
the popularity of orthogonal 
projectors in the context of approximation techniques. Nevertheless, 
suppose that the observed signal is  produced by 
the interference of 
two phenomena so that the model \eqref{eq1} becomes
\be
\label{eq2}
\sum_{i=1}^k c_i v_i  + \sum_{i=1}^n d_i w_i.
\ee
If one were interested in discriminating the 
phenomena by splitting the signal,   
the component in $\Span\{v_i\}_{i=1}^k$  could be obtained
 by an oblique projection operation mapping  the 
other component to zero. 
There is a broad range of applications in which this procedure 
happens to be of assistance \cite{sp1}. Thus,  we
felt motivated to find recursive equations 
for adapting oblique projectors. Some of the equations 
to be proposed here are inspired by our previous work on 
recursive biorthogonalization for orthogonal projectors
representation \cite{bio1,bio2}. 
We have recently been made aware that such a work is closely related 
to earlier one on recursive generalized inverses 
\cite{gre,ben,fle,moh}. 

In spite of the fact that for most numerical implementations  
a projector is
represented by a matrix, we prefer to think of projectors
as operators acting by performing inner products.
An important reason for this choice 
is the following: The equations 
can thereby be applied in general inner product spaces and 
comprise two very important cases in particular.  Namely, 
the Euclidean inner product space, where a projector is 
indeed a matrix, and the space of functions of finite 2-norm.  
We like to see  the proposed recursive equations as 
generalized Gram-Schmidt like procedures for generating sequences in 
 inner product spaces. Such sequences give rise
to oblique projectors onto nested subspaces and, of course,
to orthogonal projectors as special case.

The paper is organized as follows: Section~\ref{sec2} 
introduces the notation along with a discussion on the 
general construction of
oblique projectors. Section~\ref{sec3} provides the 
recursive equations for stepwise updating/downdating of 
such projectors. Applications are  illustrated  in 
Section~\ref{num} by i)recovering a simulated X-ray diffraction 
peak from a background and ii)filtering impulsive noise 
from the register  of the motion of a system consisting of 
the superposition of damped  harmonic oscillators. 
The conclusions are drawn in Section~\ref{con}.

\section{Oblique projectors}
\label{sec2}
As already  mentioned  we will work in a general inner product 
space $\h$, where the square norm $||.||^2$ is induced 
by the inner product that we represent as 
$\la \cdot , \cdot \ra$.  Given two closed subspaces,
$\SV \in \h$ and $\SWC \in \h$, 
such that $\h=\SV + \SWC$ and $\SV \cap \SWC = \emptyy$, 
the {\em oblique projector} operator onto $\SV$ 
along $\SWC$ will be represented as
$\EVW$. Then $\EVW$ satisfies: 
\ben
\EVW^2&=& \EVW \nonumber\\ 
\EVW v&=& v, \quad \text{for any} \quad v \in \SV  \nonumber\\
\EVW w&=& 0, \quad \text{for any} \quad w \in \SWC.  \nonumber
\een
In the particular situation in which $\SWC$ happens to be 
the orthogonal complement of $\SV$ in $\h$, i.e. if
$\SWC =\SV^\bot$, the operator is self-adjoint and represents an orthogonal 
projection onto $\SV$. We emphasize this special case by using
the particular notation $\EVV=\op_{\cal{\SV}}$. In the sequel 
the orthogonal projector operator onto a subspace,
say the subspace ${\cal{X}}$, will be indicated as
$\op_{\cal{X}}$.

Let us assume that in general $\SV=\Span\{v_i\}_{i=1}^k$ and
$\SW=\Span\{u_i\}_{i=1}^k$, with $\SW$ the orthogonal 
complement of $\SWC$.
Denoting as
$\ei, i=1,\ldots,k$ the standard orthonormal basis in 
$\C^k$, i.e., the inner product $\la \ei , {\bmth{e}_j}\ra=\delta_{i,j}$ with $\delta_{i,j}$ equal one  
if $i=j$ and zero otherwise, we 
define the operators $\hat{V}$ and $\hat{W}$
as
$$\hat{V}=\sum_{i=1}^k v_i \la \ei , \cdot \ra,
\;\;\;\;\;\;\;\;\;\;\;\;\;\;
\hat{W}=\sum_{i=1}^k  u_i \la  \ei, \cdot \ra.$$
Thus the corresponding adjoint operators  $\hat{W}^\ast$ and
$\hat{V} ^\ast$ are
$$\hat{V}^\ast=\sum_{i=1}^k \ei \la v_i , \cdot \ra,
\;\;\;\;\;\;\;\;\;\;\;\;\;\
\hat{W}^\ast=\sum_{i=1}^k \ei \la u_i , \cdot \ra.$$
The operations $\la v_i , \cdot \ra$ and 
$\la u_i , \cdot \ra$ indicate that $\hat{V} ^\ast$ and 
$\hat{W}^\ast$ act by performing  inner products in 
$\h$. The inner product is defined in such a 
way that for $f \in \h$ and $c$ a complex constant 
 the mapping $\hat{V}^\ast c f$ produces a vector in $\C^k$ of the 
form $\hat{V}^\ast c f= c \sum_{i=1}^k \ei \la v_i , f \ra$. 
The operation $\la \ei , \cdot \ra$ indicates the inner product in 
$\C^k$, thereby for ${\bmth{r}} \in \C^k$  
the mapping $\hat{V} c {\bmth{r}}$ yields a 
vector in $\SV$ of the form 
$\hat{V} c{\bmth{r}}=c\sum_{i=1}^k v_i \la \ei ,{\bmth{r}}\ra.$   
Note that the matrix representation 
of $\hat{W}^\ast \hat{V}$
has elements given by the inner products
$\la u_i,v_j\ra, i,j=1,\dots,k$. 
 The operator $$\hat{V} (\hat{W}^\ast \hat{V})^\dagger\hat{W}^\ast,$$ where $(\cdot)^\dagger$ denotes the Moore-Penrose 
 pseudo-inverse, is known to be the oblique projector 
 onto $\SV$ along $\SWC$ \cite{yon1}.
The particular choice  
\be
\label{udef}
u_i =v_i - \op_{\SWC}v_i= \op_{\SW} v_i,\quad i=1,\ldots,k
\ee
produces  the expression  
for $\EVW$ used in signal processing applications \cite{sp1}.
Certainly,  setting 
$\hat{W}= \op_{\SW}  \hat{V}$ one has the convenient equation 
\be 
\label{fo2}
\EVW=\hat{V}(\hat{V}^\ast \op_{\SW} \hat{V})^\dagger\hat{V}^\ast \op_{\SW}
\ee
that we adopt hereafter.

Amongst the many properties of oblique projectors that
have been studied we shall recall only the basic property needed for 
our purpose. It follows by applying $\op_{\SW}$ on
both sides of \eqref{fo2}, i.e.,
\be  
\label{pro}
\op_{\SW}\EVW=\op_{\SW}\hat{V}(\hat{V}^\ast \op_{\SW} \hat{V})^\dagger \hat{V}^\ast \op_{\SW}.            
\ee
Since $\hat{V}^\ast \op_{\SW}^\ast= \hat{V}^\ast \op_{\SW}$ and  $\hat{V}^\ast \op_{\SW} \hat{V}= \hat{V}^\ast \op_{\SW}^2 \hat{V}$, 
 with the substitution  $A=\op_{\SW} \hat{V}$ the right hand side of 
 \eqref{pro} takes the form $\hat{A}(\hat{A}^\ast \hat{A})^\dagger \hat{A}^\ast$. 
Such an operator is the orthogonal projector onto ${\cal{R}}(\hat{A})$. Consequently,
\be
\label{pro2}
\op_{\SW}\EVW= \op_{\SW},\quad \text{with}\quad \SW= {\cal{R}}(\op_{\SW}\hat{V})=\Span\{\op_{\SW} v_i\}_{i=1}^k.
\ee
By denoting $\ut_i= \sum_{j=1}^k {g}^\dagger_{i,j} u_j$
with ${g}^\dagger_{i,j}$ the element $(i,j)$ 
of matrix $(\hat{V}^\ast \op_{\SW} \hat{V})^\dagger $,
we can express $\EVW$ as
\be
\label{evw}
\EVW=\sum_{i=1}^k v_i \la  \ut_i , \cdot \ra.
\ee
Furthermore, from \eqref{pro2},\eqref{evw}, and \eqref{udef}  
\be
\label{euw}
\op_{\SW}= \op_{\SW}\EVW= \sum_{i=1}^k u_i \la  \ut_i , \cdot \ra.
\ee
Because $\op_{\SW}$ is self-adjoint 
$\Span\{\ut_i\}_{i=1}^k= \Span\{u_i\}_{i=1}^k=\SW$, and vice versa.
On comparing \eqref{evw} and \eqref{euw} we see that
the dual vectors $\ut_i$ are the same. This  
is of enormous assistance to derive the equations 
for adapting oblique projectors so as to 
account for the updating or downdating of the 
projecting subspace $\SV$. This will allow us to give the 
proofs of the proposed recursive equations  either 
by verification or by induction. 

\begin{remark}
It is appropriate to stress at this point that if we chose
$\SWC =\SV^\bot$ we would have
$u_i\equiv v_i,\,i=1,\ldots,k$  and consequently
$\Span\{u_i\}_{i=1}^k \equiv \Span\{v_i\}_{i=1}^k$. Hence 
 for such special situation
$\EVV\equiv \op_{\SW}\equiv \op_{\SV}$ and all the recursive equations 
of the subsequent sections would give rise to orthogonal projectors. 
\end{remark}

\section{Constructing recursive equations} 
\label{sec3}
In this section we provide the equations for 
updating and downdating oblique projectors in order to 
account for the following situations: 

Let us consider that the oblique projector $\EVKW$
onto the subspace $\SVK= \Span\{v_i\}_{i=1}^k$  along a 
given subspace $\SWC$ is known. If the subspace $\SVK$ 
is enlarged to $\SVKP$ by the inclusion of one element, i.e.,  
$\SVKP= \Span\{v_i\}_{i=1}^{k+1}$, we wish to construct 
$\EVKWP$ from the availability of $\EVKW$. On the other hand, 
if the subspace $\SVK=\Span\{v_i\}_{i=1}^k$ is reduced by the 
elimination of one element, say the $j$-th one, we 
wish to construct the corresponding oblique projector 
$\EVKWM$ from the knowledge of $\EVKW$.   
The subspace $\SWC$ is  assumed to be fixed. Its orthogonal 
complement $\SWK$ in $\hk= \SVK + \SWC$ changes with the index $k$ to
satisfy
$\hk= \SWK \oplus \SWC$, where  $\oplus$ denotes an orthogonal 
sum whilst the former is a direct sum, i.e., $\SVK\cap \SWC =
{\emptyy}$. 

\subsection{Updating the oblique projector $\EVKW$ to 
$\EVKWP$} \label{secup}
We assume that $\EVKW$ is known and write it in the explicit 
form
\be
\label{evkw}
\EVKW=\sum_{i=1}^k v_i \la  \utk_i , \cdot \ra.
\ee
Our aim is to find the vector
$\utkp_{k+1}$, and to  change the vectors 
$\utk_{i}, i=1,\ldots,k$ to 
$\utkp_{i}, i=1,\ldots,k$, so as to obtain
\be
\label{evkwp}
\EVKWP=\sum_{i=1}^{k+1} v_i \la  \utkp_i , \cdot \ra.
\ee
We will show that the duals 
$\utkp_{i}, i=1,\ldots,k+1$ can be constructed inductively 
 from  the dual of a single vector.
\begin{lemma} \label{lem1}
 For $\ut_1^1=\frac{u_1}{||u_1||^2}$, with 
 $u_1= \op_{\SW}v_1$, 
operator $v_1 \la \ut_1^1, \cdot \ra$
is the oblique projector onto the span of the 
single vector $v_1$ along $\SWC$. 
\end{lemma}
\begin{proof}
From the definition of $u_1$ (C.f. eq. \eqref{udef}) it follows that the
operator $ v_1 \la \ut_1^1, \cdot \ra = v_1 \la \frac{u_1}{||u_1||^2}, \cdot\ra $  maps  every vector
in $\SWC$ to the zero vector.  
Suppose that $f$ is in the span of $v_1$. Then $f=c v_1$ 
for some constant $c$. Since 
$\la u_1 , u_1 \ra = \la u_1 , v_1 - \op_{\SWC} v_1 \ra=
\la u_1 , v_1 \ra$ we have
$$v_1 \la \ut_1^1, c f \ra= c v_1 \frac{\la u_1, v_1 \ra}{||u_1||^2} = c v_1= f.$$
Moreover $v_1 \la \ut_1^1, v_1 \ra \la \ut_1^1,\cdot\ra 
= v_1 \la \ut_1^1, \cdot \ra$, 
which concludes the proof that $v_1 \la \ut_1^1, \cdot \ra$ is 
the oblique projector onto the span of $v_1$ along $\SWC$.
\end{proof}
In order to inductively construct
from $\ut_1^1= \frac{u_1}{||u_1||^2}$ 
the duals $\utkp_i,\,i=1,\ldots,k+1$ we have to discriminate 
 two possibilities 
\begin{itemize}
\item [i)]
$\SVKP=\Span\{v_i\}_{i=1}^{k+1}=
\Span \{v_i\}_{i=1}^{k}=\SVK$, 
i.e., $v_{k+1} \in \SVK.$
\item
[ii)]$\SVKP=\Span\{v_i\}_{i=1}^{k+1} \supset \Span \{v_i\}_{i=1}^{k}=\SVK$, i.e. $v_{k+1} \notin \SVK.$
\end{itemize}
Let us consider first the  case  i). 
Clearly if $v_{k+1} \in \SVK$ the corresponding
$u_{k+1}=v_{k+1}-\op_\SWC v_{k+1}$ 
belongs to $\SWK=\Span\{u_i\}_{i=1}^{k}$,
 because $v_{k+1}=\sum_{i=1}^k c_i v_i $ yields
$u_{k+1}=\sum_{i=1}^k c_i u_i$. 
The proposition below  prescribes how to  modify 
the corresponding dual vectors in order to guarantee 
that $\EVKWP= \EVKW$. 
\begin{proposition}
Let $v_{k+1} \in \SVK$ and vectors  $\utik$  
in \eqref{evkw} be given. For an arbitrary vector 
$y_{k+1}\in\h$ the dual vectors $\utikp$ computed as
\be
\label{utikp}
\utikp= \utik - \la u_{k+1}, \utik \ra  y_{k+1}
\ee
for $i=1,\ldots,k$ and $\ut^{k+1}_{k+1}=y_{k+1}$ 
produce the identical oblique projector as the dual vectors 
$\utik, i=1,\ldots,k$. 
\end{proposition}
\begin{proof} 

We use \eqref{utikp} to explicitly express $\EVKWP$
\ben
\label{red}
\sum_{i=1}^{k+1} v_i \la \utikp ,\cdot \ra 
&=&\sum_{i=1}^{k} v_i \la \utik , \cdot, \ra   - 
\sum_{i=1}^{k} 
v_i \la \utik , u_{k+1} \ra \la y_{k+1}, \cdot \ra +
 v_{k+1} \la  y_{k+1}, \cdot \ra\nonumber\\
&=&\EVKW - \EVKW u_{k+1} \la y_{k+1} ,\cdot\ra
+ v_{k+1} \la y_{k+1}, \cdot \ra \nonumber\\
&=&\EVKW - v_{k+1} \la y_{k+1} , \cdot \ra + 
v_{k+1} \la y_{k+1}, \cdot \ra,
\een 
where the last equality holds because $\EVKW \op_\SWC =0$ and
$\EVKW v_{k+1}= v_{k+1}$ for  $v_{k+1} \in \SVK$. 
Hence, the left hand side of \eqref{red} equals $\EVKW$.
\end{proof}
The next proposition considers the case ii)  
\begin{proposition}
Let vector $v_{k+1} \notin \SVK$ 
and vectors $\utik,i=1,\ldots,k$ in \eqref{evkw} be given. Thus
the dual vectors $\utikp$ computed as
\be
\label{eq}
\utikp=\utik - \ut_{k+1}^{k+1} \la u_{k+1}, \utik \ra, 
\ee
where $\ut_{k+1}^{k+1}= \frac{q_{k+1}}{||q_{k+1}||^2}$ with 
$q_{k+1}= u_{k+1} - \op_{\SWK} u_{k+1}$,   
provide us with the oblique projector $\EVKWP$.
\end{proposition}
\begin{proof} 
In order to organize the proof let us  establish the 
following relations: 
\ben
\label{a1}
\la q_{k+1},v_i\ra&=&0, \quad {\text{for}} \quad i=1,\ldots,k\\
\label{a2}
\la q_{k+1},v_{k+1}\ra&=&\la u_{k+1}, v_{k+1}\ra -  \la v_{k+1}, \op_{\SWC}  v_{k+1} \ra = || q_{k+1}||^2.
\een
The first relation follows from the definition of $q_{k+1}$ and the 
fact that $\op_{\SWK} v_i=u_i$ for $i=1,\ldots,k$
$$
\la q_{k+1},v_i\ra=  
\la u_{k+1},v_i\ra - \la \op_{\SWK} u_{k+1},v_i\ra=
\la u_{k+1},u_i\ra  + \la u_{k+1}, \op_{\SWC} v_i \ra - 
\la u_{k+1}, \op_{\SWK} v_i\ra= 0.
$$
On the other hand
$$
\la q_{k+1},v_{k+1}\ra=\la u_{k+1}, v_{k+1}\ra - \la u_{k+1}, \op_{\SWK} v_{k+1}\ra  = 
\la u_{k+1}, v_{k+1}\ra- \la v_{k+1}, \op_{\SWK} v_{k+1}\ra.
$$
Furthermore
\ben
||q_{k+1}||^2&=& \la q_{k+1},u_{k+1} \ra - \la q_{k+1},\op_{\SWK} u_{k+1} \ra =
\la q_{k+1} , u_{k+1} \ra \nonumber \\
&=&\la q_{k+1} , v_{k+1} \ra -  \la q_{k+1}, \op_{\SWC} v_{k+1}\ra =   
\la q_{k+1}, v_{k+1}  \ra \nonumber\\
&=& \la u_{k+1}, v_{k+1} \ra  -  \la v_{k+1}, \op_{\SWK} v_{k+1} \ra. \nonumber 
\een
We are now in a position to start the proof of the proposition by 
induction. From Lemma \ref{lem1} we know that
$v_1 \la u^1_1 , \cdot \ra/||u^1_1||^2 $ is the oblique
projector onto $\SV_1$ along $\SWC$. Assuming that
$\EVKW$ is the oblique projector onto $\SVK$ along $\SWC$
we will prove that $\EVKWP$ is the oblique projector
onto $\SVKP$ along $\SWC$. For this we need to prove that 
 the recursive equation \eqref{eq} yields the operator 
 $\EVKWP$ satisfying:
 \begin{itemize}
 \item [i)] 
 $\EVKWP^2= \EVKWP$ 
 \item [ii)]
 $\EVKWP v= v, \quad \text{for any} \quad v \in \SVKP $ 
 \item [iii)]
 $\EVKWP w= 0, \quad \text{for any} \quad w \in \SWC. $ 
 \end{itemize}
 We begin by using \eqref{eq}
 to express $\EVKWP$ as
 \ben
 \label{obkp}
 \sumikp v_i \la \utikp , \cdot \ra &=&
 \sumik v_i \la \utik , \cdot, \ra   -
 \sumik v_i \la \utik , u_{k+1} \ra \la \ut_{k+1}^{k+1},
 \cdot \ra + v_{k+1} \la  \ut_{k+1}^{k+1}, \cdot \ra \nonumber\\
 &=& \EVKW - \EVKW u_{k+1}\la \frac{q_{k+1}}{||q_{k+1}||^2},\cdot\ra + v_{k+1} \la \frac{q_{k+1}}{||q_{k+1}||^2},\cdot\ra.
 \een
For all $w$ in $\SWC$ it holds that $\EVKW w =0$ and
$\la q_{k+1}, w \ra=0$. Then from \eqref{obkp} 
we conclude that 
 condition iii) is satisfied. 
Every $v \in \SVKP$ can be written 
as $v=\sum_{i=1}^{k+1} c_i v_i$. Thus, 
from \eqref{obkp} and using relations \eqref{a1} and \eqref{a2}
\ben
\EVKWP v &=& \sum_{i=1}^{k} c_i v_i + c_{k+1}\EVKW v_{k+1} - c_{k+1} \EVKW u_{k+1} + c_{k+1} v_{k+1} \nonumber\\
       &=&\sum_{i=1}^{k} c_{i} v_i  + c_{k+1} v_{k+1}= v,\nonumber
\een
which demonstrates condition  ii). Finally, 
since from \eqref{obkp} and \eqref{a1} it 
follows that $\EVKWP \EVKW = \EVKW$, we have 
\ben
\EVKWP^2 &=& \EVKW - 
\EVKW u_{k+1} \la\frac{q_{k+1}}{||q_{k+1}||^2},\cdot\ra
+ \EVKWP v_{k+1}
 \la\frac{q_{k+1}}{||q_{k+1}||^2},\cdot\ra\nonumber\\
 &=& \EVKW - \EVKW u_{k+1} \la \frac{q_{k+1}}{||q_{k+1}
 ||^2},
 \cdot\ra +v_{k+1} \la\frac{q_{k+1}}{||q_{k+1}||^2},\cdot\ra
 = \EVKWP.
\een
\end{proof}

\begin{property}
\label{re1}
If vectors $\{v_i\}_{i=1}^k$ are linearly independent 
they are also biorthogonal to the dual 
vectors arising inductively from 
the recursive equation \eqref{eq}.
\end{property}
The proof of 
this  property  is  given in  Appendix A.
\begin{remark}
\label{re2}
If vectors $\{v_i\}_{i=1}^k$ are not linearly independent 
the oblique projector $\EVKW$ is not unique. Indeed, if
 $\{\utik\}_{i=1}^k$ are dual vectors giving rise to $\EVKW$ 
then one can construct infinitely many  duals as:
\be
\label{nonu}
\til{y}_i= \utik + y_i - \sum_{j=1}^k y_j \la v_j , \utik \ra
\quad i=1,\ldots,k,
\ee
where $y_i,\,  i=1,\ldots,k$  are arbitrary vectors 
in $\h$.
\end{remark}
\begin{proof}
We use \eqref{nonu} to write
\ben
\sumik v_i \la  \til{y}_i , \cdot \ra & =&
\EVKW + \sumik v_i \la {y}_i , \cdot \ra - 
\sumik v_i \sum_{j=1}^k  \la \utik,  v_j \ra \la  y_j, \cdot\ra  
\nonumber\\
&=& \EVKW + \sumik v_i \la {y}_i , \cdot \ra -
\sum_{j=1}^k\sumik  v_i \la  \utik,  v_j \ra \la {y}_j 
, \cdot \ra \nonumber \\
&=&\EVKW + \sumik v_i \la {y}_i , \cdot \ra - \sum_{j=1}^k
 v_j \la {y}_j , \cdot \ra \nonumber \\
& =&\EVKW.
\een
\end{proof}
It follows from Property \ref{re1} that if 
vectors $\{v_i\}_{i=1}^k$ are linearly independent 
equation \eqref{nonu} yields the unique duals
 $\til{y}_i \equiv \utik,\, i=1,\ldots,k$.
\subsection{Downdating the oblique projector $\EVKW$ to
$\EVKWM$}
Suppose that by the elimination of 
the element $j$ the subspace $\SVK$ 
is reduced to $\SV_{\kj}= \Span\{v_i\}_{\inoj}^k$. In order 
to give the equations for adapting the corresponding dual 
vectors generating  the oblique projector $\EVKWM$ we 
need to consider two situations: 
\begin{itemize}
\item
[i)]$\SVKM=\Span\{v_i\}_{\inoj}^{k}=\Span\{v_i\}_{i=1}^{k}=\SVK$
i.e., $v_j \in \SVKM.$
\item
[ii)]$\SVKM=\Span\{v_i\}_{\inoj}^{k} \subset \Span\{v_i\}_{i=1}^{k}=\SVK$, i.e., $v_j \not\in \SVKM.$
\end{itemize}
The next proposition addresses i).
\begin{proposition}
Let $\EVKW$ be given by \eqref{evkw} and let us assume 
that removing vector $v_j$ from the spanning set of 
$\SVK$ leaves the identical subspace, i.e., $v_j \in \SVKM$.
Hence, if the remaining dual vectors are 
modified as follows: 
\be
\label{opeqb} 
\utikj= \utik + 
\frac{\la u_j, \utik\ra  \ut_j^k}
{1- \la u_j, \ut_j^k\ra},
\ee
the corresponding oblique projector does not change, 
i.e.  $\EVKWM=\EVKW$.
\end{proposition}
\begin{proof} 
Let us recall that $v_j \in \SVKM$ implies
$u_j \in \SVKM$. 
 Hence $\la u_j, \ut_j^k\ra \neq 1$, as it is seen 
from the fact that
$\la u_j, \ut_j^k\ra = \sum_{i=1}^k \la u_j, \ut_i^k \ra
\la u_i, \ut_j^k \ra= \sum_{i=1}^k \la u_j, \ut_i^k \ra
 \la \ut_i^k, u_j \ra=
\sum_{i=1}^k | \la u_j, \ut_i^k\ra|^2$, which implies
$\la  u_j, \ut_j^k\ra=1$ if and only if for
$i=1,\ldots,k$ it holds that $\la u_j, \ut_i^k\ra=
\delta_{i,j}$. This is not true if $v_j \in \SVKM$ 
so that we can use \eqref{opeqb}  to 
express $\EVKWM$ as 
\ben
\sum_{\inoj}^k v_i \la \utikj, \cdot \ra & = &
\sum_{\inoj}^k v_i \la \utik, \cdot\ra + 
\sum_{\inoj}^k  
\frac{ v_i \la  \utik,  u_j\ra \la \ut_j^k,\cdot\ra }
{1- \la \ut_j^k, u_j\ra}\nonumber\\
&=&\EVKW - v_j \la \ut_j^k,\cdot \ra +
 \frac{ \EVKW u_j\la \ut_j^k,\cdot\ra}{1- \la \ut_j^k,u_j\ra} -  v_j \la \ut_j^k, u_j  \ra   \frac{\la \ut_j^k,\cdot\ra}
{1- \la   \ut_j^k, u_j \ra} \nonumber\\
&=&\EVKW - v_j \la \ut_j^k,\cdot \ra + \frac{v_j\la \ut_j^k, \cdot\ra} {1- \la \ut_j^k, u_j\ra} 
-  \frac{v_j \la \ut_j^k,\cdot\ra }{1- \la \ut_j^k, u_j\ra} \la  \ut_j^k ,u_j  \ra  \nonumber\\
&=& \EVKW. 
\een
\end{proof}
Finally, proposition 4 addresses ii).
\begin{proposition}
Let $\EVKW$ be given by \eqref{evkw} and let us assume
that the vector $v_j$ to be removed from the spanning 
set of $\SVK$ is not in $\SV_{\kj}$. 
In order to produce the oblique projector $\EVKWM$ the 
appropriate modification of the dual vectors  can be achieved 
by means of the following equation
\be
\label{duba}
\utikj=\utik-\frac{\ut_j^k \la \ut_j^k , \utik\ra}{||\ut_j^k||^2}.
\ee
\end{proposition}
\begin{proof} Using \eqref{duba} we write:
\ben \label{equp}
\sum_{\inoj}^k v_i \la \utikj, \cdot \ra &=& 
\sum_{\inoj}^k v_i \la \utik , \cdot \ra - 
\sum_{\inoj}^k  \frac{v_i \la \utik , \ut_j^k \ra 
\la \ut_j^k, \cdot \ra}{||\ut_j^k||^2} \nonumber\\ 
&=&\EVKW - v_j \la \ut_j^k , \cdot \ra -\frac{ \EVKW \ut_j^k \la \ut_j^k, \cdot \ra}{||\ut_j^k||^2} + 
v_j \la \ut_j^k , \cdot \ra \nonumber\\
&=&\EVKW - \frac{ \EVKW \ut_j^k \la \ut_j^k, \cdot \ra}{||\ut_j^k||^2}. 
\een
We notice that $\frac{\ut_j^k \la \ut_j^k, \cdot \ra}{||\ut_j^k||^2}$ is the 
orthogonal projector onto the span of the single vector
$\ut_j^k$ and denote it as $\op_{\utk_j} = \frac{\ut_j^k \la \ut_j^k, \cdot \ra}{||\ut_j^k||^2}$.
Thus the orthogonal projector onto ${\SW_{\kj}}$ can be  expressed as
$\op_{\SW_{\kj}}= \op_{\SWK} - \op_{\utk_j}$. 
Applying this operator on the right hand side of \eqref{equp} 
we obtain:
$$\op_{\SW_{\kj}} (\EVKW-\EVKW \op_{\ut_j^k})= (\op_{\SWK} - \op_{\utk_j}) (\EVKW-\EVKW \op_{\ut_j^k})$$
and, 
since $\op_{\SWK} \EVKW= \op_{\SWK}$ and $\op_{\utk_j} \op_{\SWK} = \op_{\utk_j}$,  using 
the fact that $ \op_{\utk_j} \EVKW= \op_{\utk_j} \op_{\SWK} \EVKW= \op_{\utk_j}$ it follows that
\ben 
(\op_{\SWK} - \op_{\utk_j}) (\EVKW-\EVKW \op_{\utk_j})& =& 
\op_{\SWK} - \op_{\utk_j} -  \op_{\utk_j} + 
\op_{\utk_j} \nonumber\\
& =&\op_{\SWK} - \op_{\utk_j} = \op_{\SW_{\kj}} = \op_{\SW_{\kj}} \EVKWM. \nonumber
\een
From the last equation and \eqref{equp} we gather that
\be
\op_{\SW_{\kj}} \sum_{\inoj}^k v_i \la \utikj, \cdot \ra   
- \op_{\SW_{\kj}}  \EVKWM =0.  
\ee
For every vector $f\in \h$ we  therefore have
\be
\op_{\SW_{\kj}}(\sum_{\inoj}^k v_i \la \utikj, \cdot \ra
- \EVKWM) f = \op_{{\SW}_{\kj}} \Delta \hat{D} f =0,
\ee
with $\Delta \hat{D}= \sum_{\inoj}^k v_i \la \utikj, \cdot \ra - \EVKWM$. 
This implies that either $\Delta \hat{D}$ is the zero operator or 
$\Delta \hat{D} f \in \SWC$ 
for every $f \in \h$. The latter cannot be true because from the
definition of $\Delta \hat{D}$ it is seen that $\Delta \hat{D} f \in \SVKM$ 
and by hypothesis $\SVKM \cap \SWC =\emptyy$. 
Hence  $\Delta \hat{D}$ should be the zero operator, 
which leads to  the conclusion that
\be
\EVKWM=\sum_{\inoj}^k v_i \la \utikj, \cdot \ra.
\ee
\end{proof}
\begin{remark}
\label{re3} 
The case of replacing a vector in $\SV_k$, say $v_j$ by $v'_j$,  
is actually equivalent to augmenting the 
subspace $\SVKM$ to $\SVKM  + v'_j$ after the vector 
$v_j$ was deleted. Some implementation 
issues arise, though. In order to modify the duals as 
prescribed in \eqref{eq} we need to compute a vector 
$q_{\kj}=u'_j - \op_{\SW_{\kj}} u'_j$. For the sequential enlargement 
of the projecting subspace, discussed in section \ref{secup}, 
the projector $\op_{\SW_k}$ can be sequentially constructed 
by means of the orthonormal vectors 
$q_{n}/||q_n||,\, n=1,\ldots,k$.
Nevertheless, when replacing vectors sequentially we need to 
allow for the recalculation of the corresponding projectors. 
One possibility that could be considered is 
the recalculation of the orthogonal vectors $q_{n}$ \cite{rei,boj}.
An alternative approach implies to use of the dual 
corresponding to the deleted vector for orthogonalization purposes. 
A discussion concerning the implementation of such a procedure 
is given in \cite{andre2,andre3}. 
\end{remark}

\section{Applications to signals discrimination}
\label{num}
The examples presented in this section aim at illustrating
 the application of our recursive construction of 
oblique projectors for signals in $L_2[a\,,\,b]$, 
the space of square integrable functions on $[a\,,\,b]$. 
For $f$ and $g$ in  $L_2[a\,,\,b]$,  
we define the inner product, according to the previously 
adopted convention, as
$$\la f , g \ra=\int_{a}^ {b} f^\ast(x) 
g(x)\,dx,$$
where $f^\ast(x)$ indicates the complex conjugate of $f$. 
In the examples bellow  all the integrals are numerically 
calculated.

\subsection{Extraction of a X-ray diffraction peak from
a dispersive background}

Here the signal is simulated by emulating a crystallographic 
problem. 
It is assumed to be the X-ray diffraction
intensity produced by a powder sample of a clay mineral
consisting of very flat crystals.  Each  such crystal is formed
by the stacking of $n$ layers  producing a diffraction
intensity as given by \cite{klu,jen}
\be
I_n(x) = \frac{\sin^2 n x }{\sin^2 x},
\ee
where the variable $x$, given in radians,
is related to the diffraction
angle $\theta$ according to the equation
\be
x= 2\pi \frac{d}{\lambda} \sin  \theta.
\label{xt}
\ee
The parameter $d$ in \eqref{xt} is the effective distance between two consecutive 
layers and  characterizes the material. The parameter    
$\lambda$ represents the wavelength  of the incident
radiation. 

We denote the diffraction intensity produced by the whole sample 
as $f_1(x)$. Thus, 
\be
\label{smod}
f_1(x)= \sum_{n=1}^{k} c_n \frac{\sin^2 n x }{\sin^2 x}.
\ee
As already mentioned $n$ indicates the possible number
of layers forming
a single crystal in the sample. The coefficients $c_n$
account for the proportions of crystals consisting of $n$ layers.
Here, for simulating the signal,  the coefficients were
considered to be
$$c_n= e^{-0.05(n-7)^2} + 0.2 e^{-0.1(n-35)^2},\,
n=1\ldots,60.$$
The diffraction figure  $f_1$ emerges from a background
that is modelled as
\be
f_2(x)=  50 \sum_{j=1}^{3} j e^{-j(x-\frac{\pi}{2})}.
\ee
The combined phenomenon gives rise to the signal
$f=f_1+f_2$ plotted in the left graph of Figure~1 on the 
interval relevant to the diffraction  model, namely 
$[\frac{\pi}{2}\,,\, \frac{3\pi}{2}].$

We are interested in extracting the diffraction peak $f_1$ 
from the background. For this we will construct sequentially 
oblique projectors onto subspaces $\SVK$ 
given as
$$\SVK=\Span \{ \frac{\sin^2 n x }{\sin^2 x},\, n=1,\ldots,k\},\,
x \in [\frac{\pi}{2}\,,\, \frac{3\pi}{2}].$$
The final $k$-value is to be adjusted. 
The subspace $\SWC$ is here
$$\SWC= \Span \{e^{-j(x-\frac{\pi}{2})}, \, j=1,\ldots,3\},\, 
x \in [\frac{\pi}{2}\,,\, \frac{3\pi}{2}].$$ 
Since the order $k$ of the diffraction model is assumed unknown, 
it was adjusted as follows: firstly the order model 
was sequentially  increased 
(up  to $k=200$) and then sequentially downdated. 
It was observed 
that the recovering of the signal was not very sensitive 
to the model order. In a range from $k=50$ to $k=200$ 
the approximations were  totally equivalent. 
From $k=40$ to $k=50$ changes in the approximations 
 were noticed  but the approximations
could still be considered `practically' equivalent. 
The recovered peak $f_1= \EVKW f$ (for $k=50$) is depicted 
in the right graph of  Figure~1. It happens to coincide, 
in the scale of the 
figure, with the graph of the theoretical one (C.f. \eqref{smod}). 

The convenience of the proposed 
adaptive technique in the determination of the order 
of the diffraction intensity model is clear:  otherwise when changing 
the $k$ value as described above the whole projector
would have to be recalculated for each different value of $k$. 
However, the advantage of the proposed technique is even more significant 
when, as is the case in the next example, overestimation 
of the order in the signal model may result in the failure to 
discriminate the signals.


\begin{figure}[!ht]
\label{f1}
\begin{center}
\includegraphics[width=7.5cm]{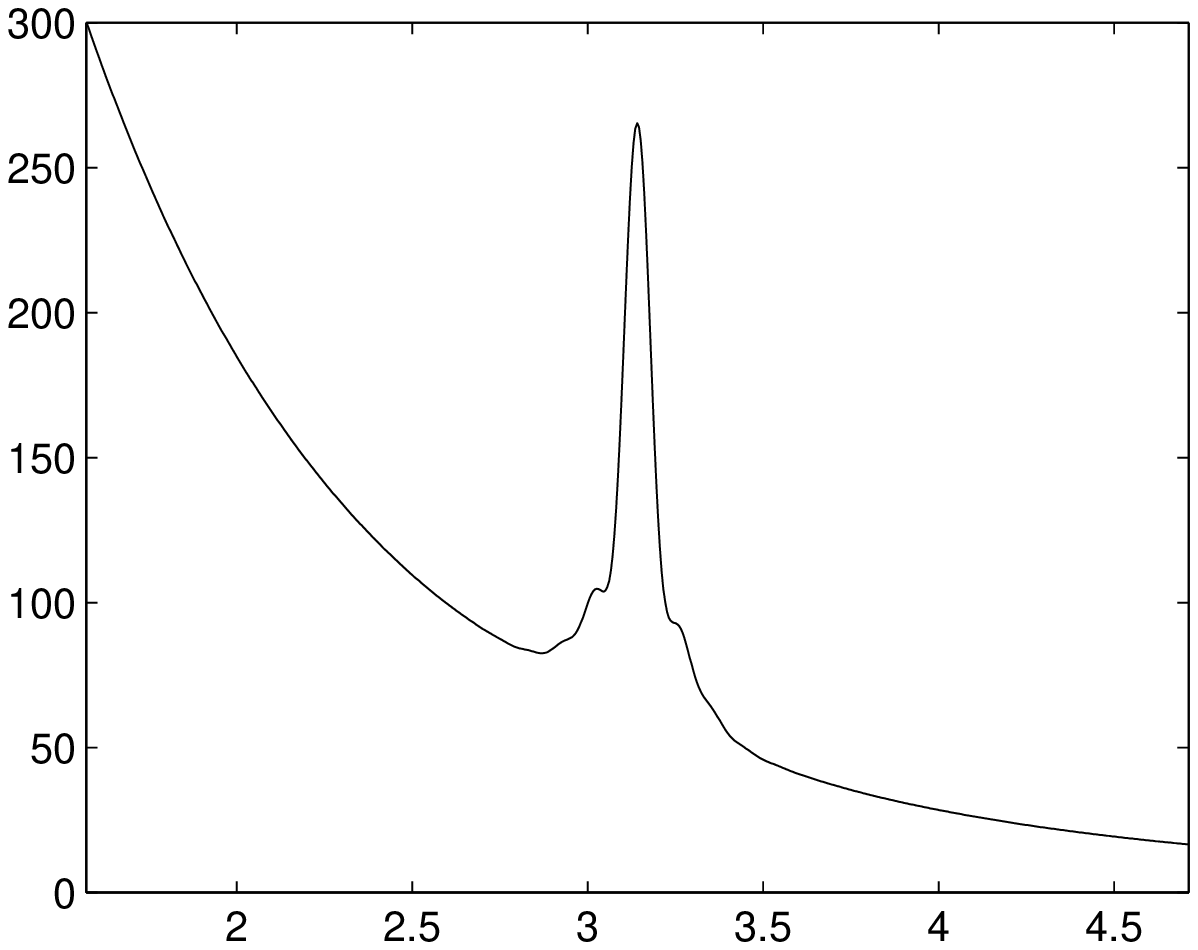}   
\includegraphics[width=7.5cm]{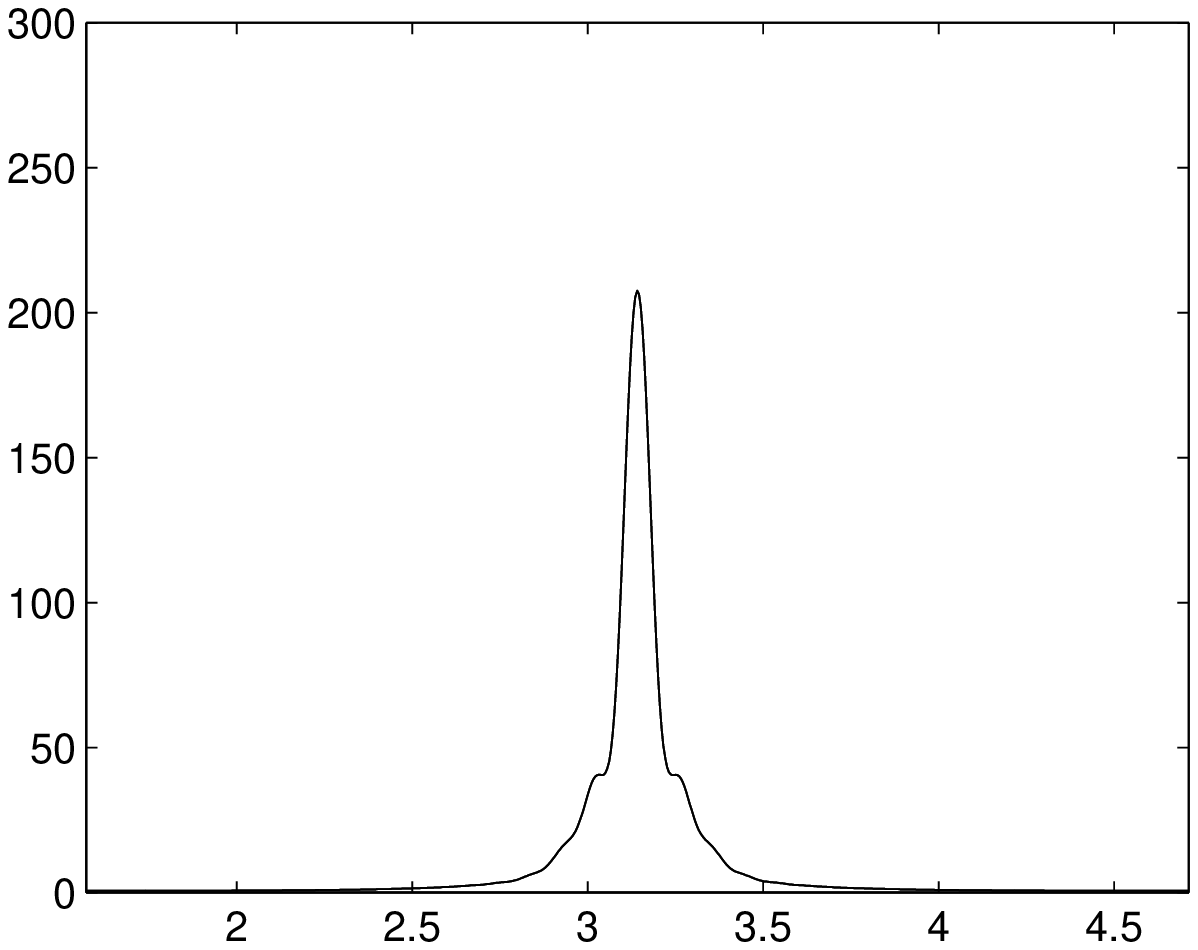}
\end{center}
\caption{Graphs on the left: Diffraction intensity 
vs the variable $x$ measured in 
radians and related to the diffraction angle through \eqref{xt}.
Graph on the right: Diffraction intensity extracted from 
the background.}
\end{figure} 

\subsection{Elimination of impulsive noise}
In this case the signal is considered to be the register
of the motion of a system consisting of 
uncoupled damped harmonic oscillators.
The $n$-th oscillator is characterized by a
frequency of $\frac{n}{2}$ Hz. The corresponding
equation for its motion $x_n(t)$ as a function of time 
is given  as
\be
x_n(t)= e^{-t} \cos (\pi n t), \quad n=1,\ldots,k.
\ee
The distribution of frequencies is considered to be
$c_n= (1+0.7(n-75)^2)^{-1}$ so that the motion of the system is
registered by the signal
\be \label{smod2}
f_1(t)=\sum_{n=1}^{100} 
\frac{e^{-t} \cos (\pi n t)} {1+0.7(n-75)^2},\quad t \in [0,1].
\ee 
This signal (shown in  the right 
graphs of Figure 2) is corrupted by impulsive noise, 
which represents a type of electrical noise appearing in 
some practical situations. 
The possible pulses are taking from the set of $400$
 Gaussian sparks  
$p_j(t)= e^{-100000(t-0.0025j)^2},\, j=1, \ldots,400$. Hence the 
corresponding subspace $\SWC$ is 
$$\SWC= \Span \{e^{-100000(t-0.0025j)^2},\, j= 1,\ldots,400\},\quad
 t \in [0,1].$$ 
First a random superposition of 
$50$ pulses is added to the signal $f_1$ to simulate the 
noisy one plotted in the top 
left graph of Figure 2. The result after the oblique projections 
along the subspace $\SWC$ given above is shown in the 
top right graph (it coincides with the theoretical signal $f_1$ 
given in \eqref{smod2}). 
The left graph at the bottom of Figure 2 corresponds  to 
a different realization of the noise, in this case 
generated as a random superposition of $200$ pulses. 
The signal, after filtering  by  oblique projections
along $\SWC$, is shown in the right graph (it also coincides with 
the theoretical one).  
\begin{figure}[!ht]
\label{f2}
\begin{center}
\includegraphics[width=7.5cm]{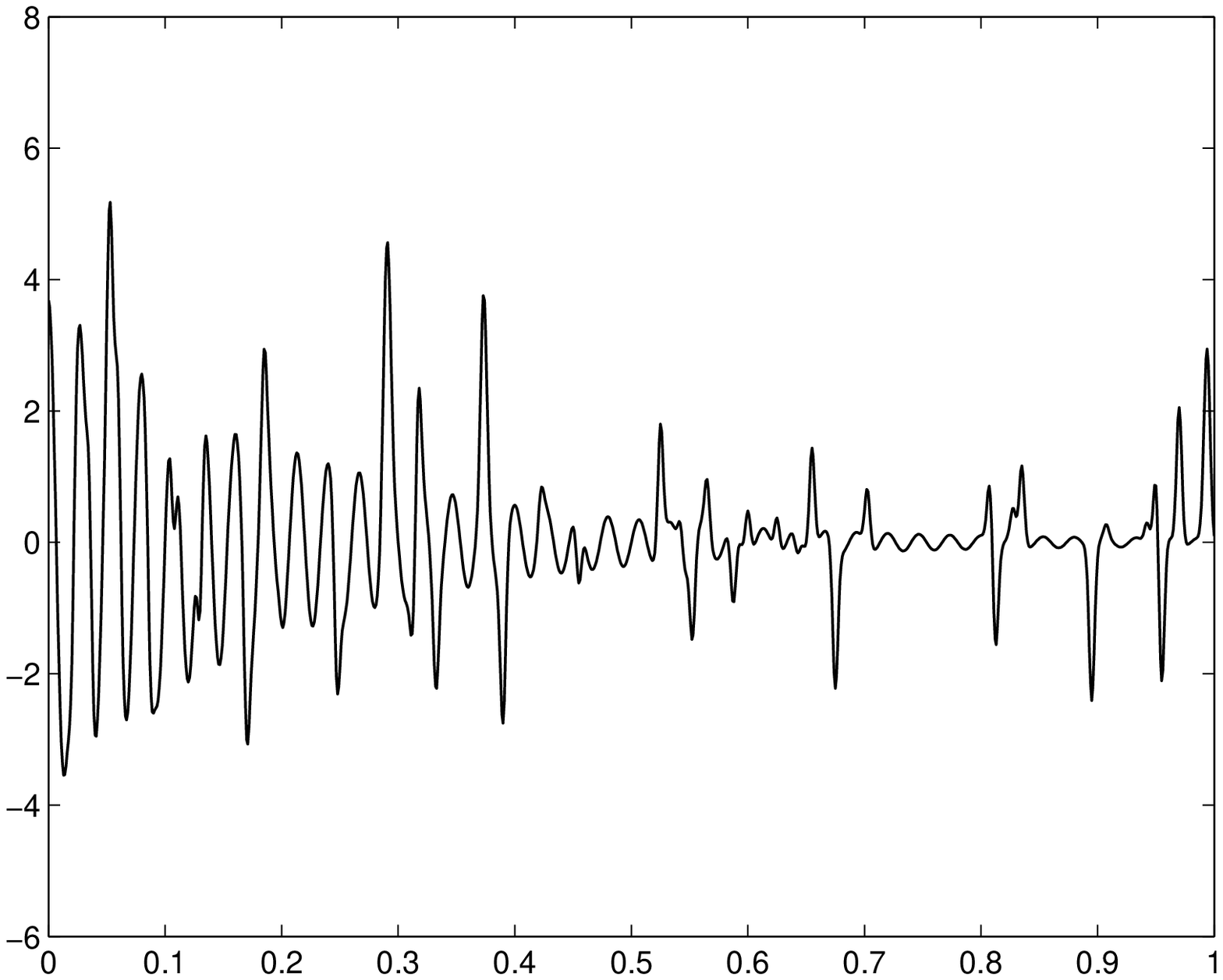}
\includegraphics[width=7.5cm]{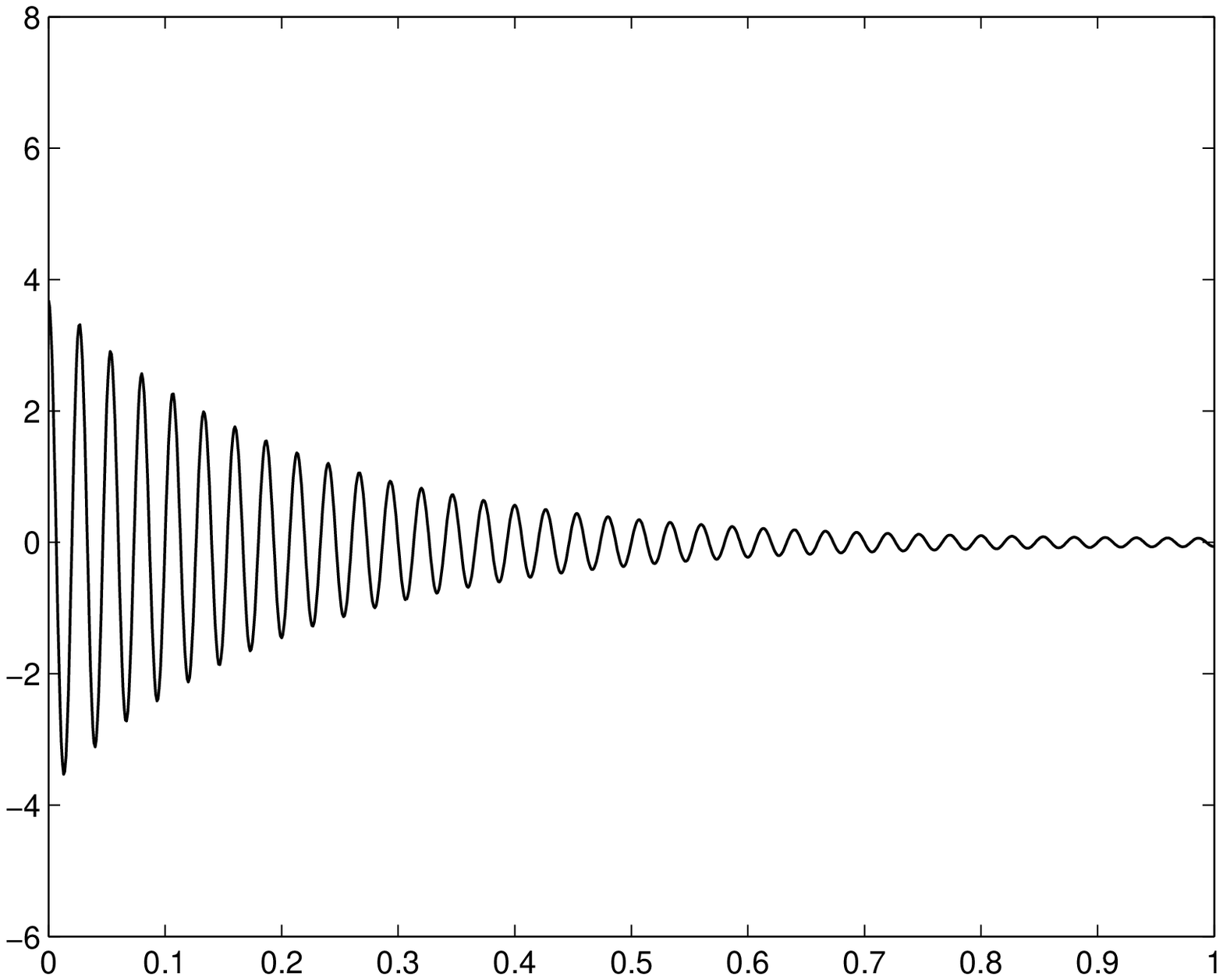}\\
\includegraphics[width=7.5cm]{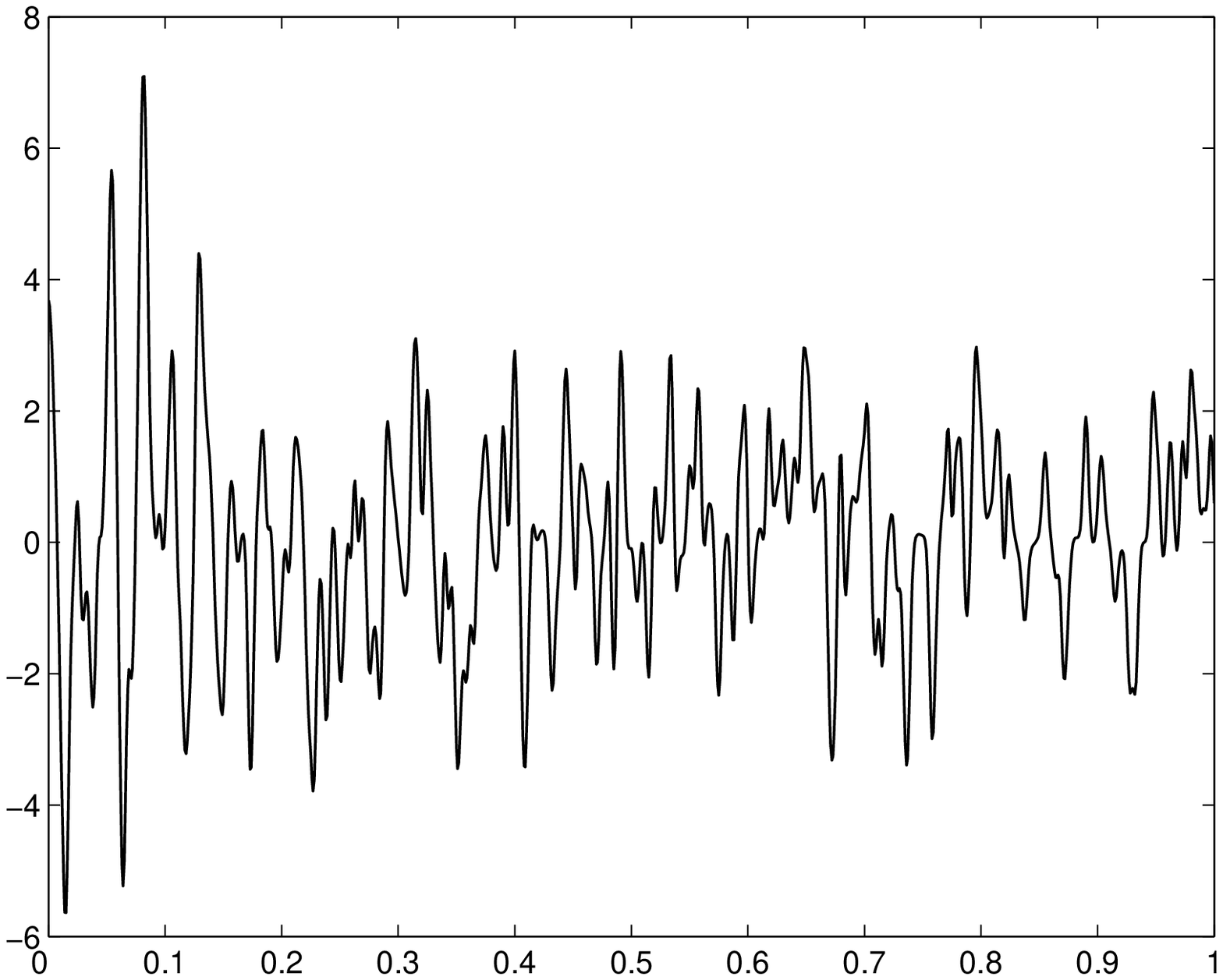}
\includegraphics[width=7.5cm]{osc.eps}\\
\end{center}
\caption{Top left graph: motion of the harmonic oscillators system 
as a function of time (in seconds) corrupted by 50 random pulses. 
The graph on the right depicts the signal after filtering the noise 
by sequential oblique projection. 
The bottom figures have the same description but the noise corresponds 
to 200 random pulses.} 
\end{figure}

Let us point out that an alternative
way of splitting the signal would entail fitting both the
signal and noise models, with the consequent increment in
the dimension of the problem of determining the
corresponding unknown parameters.
In the example of this section,
for instance,  400 more parameters (coefficients of the noise model)
would be involved.

\section{Conclusions}
\label{con}
Recursive equations for updating/downdating
oblique projectors have been proposed. The updating
strategy can be regarded as a generalized
Gram-Schmidt like  procedure for generating
a sequence giving rise to oblique projectors along
a fixed subspace. The downdating strategy modifies 
such sequence to account for the removal of some 
elements.
The equations are of the same nature as those
for producing orthogonal projectors,
but involve different vectors. Orthogonal
projectors arise within this framework  as
a particular case. 

The proposed technique has  been applied to the 
problem of discriminating signals produced by 
different phenomena. The applications that have been considered 
assume that the signal model is  determined by 
physical considerations and only the order of the model 
is to be adjusted. The task of setting the right order model 
is facilitated by the recursive nature of the proposed
equations. In the two examples considered here the signal 
splitting is not very sensitive to the order of  
the signal model. However, an important difference 
between the two examples is the following: 
while in the first example an excessive overestimation of order model
(maximum possible number of layers present in a crystal) 
does not prevent the extraction of the diffraction peak
from the background, an excessive overestimation of the 
order of the  model in the second example  
(maximum possible frequency of an oscillator) may produce the  
failure to separate the signal from the impulsive noise. 
The reason being that  
for very high frequencies the angle between the signal subspace
and the noise subspace becomes very small, which generates an 
ill posed problem. 

The recursive feature of the proposed equations turns out to be   
even more important in those situations 
in which the signals splitting is achieved by stepwise selection of 
each component of the signal model. 
This is the subject of a recent work
\cite{oblmp}, where the present approach is  
shown to be of significant assistance. 

\subsection*{Acknowledgements}
Support from EPSRC (EP$/$D062632$/$1) is acknowledged.
\newpage

\appendix
\section*{Appendix A: Proof of Property \ref{re1}}
\renewcommand{\theequation}{A.\arabic{equation}}
\setcounter{equation}{0}
\label{apen}
Let us recall that if vectors $\{v_i\}_{i=1}^{k+1}$ are linearly 
independent, all the duals $\{\ut^{k+1}_i\}_{i=1}^{k+1}$ are generated by 
the recursive equation \eqref{eq}. We need to show that such 
vectors satisfy:
$$
\la v_m, \utkp_i \ra = \delta_{m,i},\quad m,i=1,\ldots,k+1.
$$
\begin{proof}
For $k=0$ the relation holds because 
$\ut^1_1 = \frac{u_1}{||u_1||^2}$ and 
$||u_1||^2 = \la v_1, v_1 \ra - \la v_1, \op_{\SWC} v_1 \ra$. 
Therefore 
$\la v_1, \ut^1_1 \ra= 
\frac{\la v_1, {v_1}\ra - \la v_1, \op_{\SWC} v_1 \ra} 
{\la v_1, v_1 \ra - \la v_1, \op_{\SWC} v_1 \ra} = 1$. 

Assuming that for 
$k+1=l$ it is true that 
$$
\la v_m, \ut^l_i \ra = \delta_{m,i},\quad m,i=1,\ldots,l
$$
we will prove that
$$
\la v_m, \ut^{l+1}_i \ra = \delta_{m,i},\quad m,i=1,\ldots,l+1.
$$
For this we need to consider four different situations 
with regard to the indices. 
\begin{itemize}
\item [I)] $m=1,\ldots,l$ and $i=1,\ldots,l.$

In this case 
$\la v_m, q_{l+1} \ra = 0$ (C.f. \eqref{a1}).
Hence, from the recursive equation
\eqref{eq}, we have
$$\la v_m, \ut^{l+1}_i\ra = \la v_m, \ut^{l}_i \ra +0 
=\delta_{m,i}.$$
\item [II)] $m=l+1$ and $i=1,\ldots,l$

Now
\ben
\la v_{l+1}, \ut^{l+1}_i \ra &= &
\la v_{l+1}, \ut^{l}_i \ra  - \frac{\la v_{l+1}, q_{l+1}\ra}
{||q_{l+1}||^2} \la v_{l+1}, \ut^{l}_i \ra \nonumber\\
&= & \la v_{l+1}, \ut^{l}_i \ra \frac{||q_{k+1}||^2 - 
\la v_{l+1}, q_{l+1}\ra} {||q_{l+1}||^2}\nonumber
\een
so that, since $\la v_{l+1}, q_{l+1}\ra= ||q_{k+1}||^2$ 
(C.f. \eqref{a2}), 
$$\la v_{l+1}, \ut^{l+1}_i \ra = 0.$$
\item[III)]$m=l+1$ and $i=l+1.$

This implies
$$\la v_{l+1}, \ut^{l+1}_{l+1} \ra= 
\frac{\la v_{l+1}, q_{l+1}\ra}{||q_{l+1}||^2}=1.$$

\item[IV)]$m=1,\ldots,l$ and $i=1+l.$ 

In this case
$$\la v_{m}, \ut^{l+1}_{l+1}\ra =
\frac{\la v_{m}, q_{l+1}\ra}{||q_{l+1}||^2}=0.$$ 

\end{itemize}
From I) II) III) and VI) we conclude that  
$$\la v_m, \ut^{l+1}_i \ra = \delta_{m,i},\quad m,i=1,\ldots,l+1,$$
which proves that the vectors generated through \eqref{eq}  are
 biorthogonal to vectors $v_{m}, m=1,\ldots,k+1$. 
\end{proof}

\end{document}